\documentclass[11pt]{article}
\usepackage{amsfonts,amssymb}
\usepackage{epsfig}
\usepackage{latexsym}
\usepackage[all]{xy}
\usepackage{amsmath}
\usepackage{amsfonts}
\usepackage{amssymb}
\usepackage{fullpage}

\def\R{\mathbb R}

\def\Z{\mathbb Z}

\def\S{\mathbb S}
\def\C{\mathbb C}

\def\Q{\mathbb Q}

\def\O{\mathbb O}
\def\H{\mathbb H}

\newtheorem{Thm}{Theorem}[section]
\newtheorem{Cor}{Corollary}[section]
\newtheorem{Lemma}{Lemma}
\newtheorem{Prop}{Proposition}[section]

\newtheorem{Remark}{Remark}[section]
\begin{document}
\title{Combinatorial 8-Manifolds having Cohomology of the Quaternionic Projective Plane and Their Nonembeddings}
\author{Satya Deo\thanks{This work was done while the author was supported by the DST (Govt of India) research grant sanction letter number SR/S4/MS:567/09 dated 18.02.2010}\\Harish-Chandra Research Institute \\
Chhatnag Road, Jhusi,\\
Allahabad 211 019, India.\\
E-mail: {\tt sdeo@mri.ernet.in,  vcsdeo@yahoo.com}\\[2mm] }
\date{}
\maketitle

\begin{abstract}
In this paper we prove, using $\Z_2-$index theory, that none of the three combinatorial 8-manifolds on 15 vertices constructed by Brehm and  K\"{u}hnel, each of which is a cohomology quaternionic projective plane, can be combinatorially embedded in the Euclidean space $E^{12}$, though they have tight polyhedral embeddings in $E^{14}.$ This extends a similar method already known for the nonembeddings of real and complex projective planes.
\end{abstract}

\section{Introduction}
It is well known that the real projective plane $\R P^2$ cannot be embedded in the 3-dimensional Euclidean space $E^3$, though it has a tight polyhedral embedding in $E^5$. Similarly, the complex projective plane $\C P^2$ cannot be embedded in the 6-dimensional Euclidean space $E^6$~\cite{mat}, though it has a tight polyhedral embedding in $E^8$~\cite{mat}. In view of these results, it is natural to investigate the nonnembedding questions for the quaternionic projective plane $\H P^2$, and also for the Cayley plane (Octonionic Projective plane) $\O P^2$ . Using the characteristic class theory for embeddings, it is, of course, known that neither $\H P^2$ nor $\O P^2$ can be topologically embedded in $E^{12}$ and $E^{24}$ respectively, but the proofs use the advanced machinery of characteristic classes in algebraic topology. Here we use the idea of $\Z_2-$index theory to give an easy geometric proof that none of the three combinatorial manifolds, called BK complexes, which are cohomology quaternionic projective plane, can be combinatorially embedded in $E^{12}$.

We know that $\H P^2$ is defined as the quotient space of the 2-dimensional quaternionic sphere $\{(q_0,q_1,q_2)\;|\; ||q_0||^2+||q_1||^2+||q_2||^2 =1\;\hbox{where $q_0,q_1,q_2$ are quaternionic numbers}\}$ under the action of unit quaternionic sphere $\S^3=\{q\in\Q|\;||q||=1\}$ acting freely on right, i.e., $\H P^2=\S^{11}/\S^3$. Therefore, it follows from the definition that $\H P^2$ is an 8-dimensional, compact, simply connected real manifold. It is well known that $\H P^2$ is a CW-complex having one $i$-cell for $i=0,4,8$ and that its cohomology algebra over $\Z$ is the truncated polynomial algebra $$H^* (\H P^2,\Z)=\Z[a]/(a^3)$$ where the degree of generator $a$ is 4 (see~\cite{hat}). It is also known that $\H P^2$ admits a finite simplicial complex structure, and it has been proved by Brehm and K\"{u}hnel~\cite{bk1} that any triangulation of $\H P^2$ must have at least 15 vertices. Furthermore, Brehm and K\"{u}hnel  have also explicitly constructed three 8-dimensional combinatorial manifolds $M^{8}_{15},\;\widetilde{M}^{8}_{15},\;\widetilde{\widetilde M^{8}_{15}}$ which have the same cohomology algebra as $\H P^2$. We will call these simplicial complexes as {\it BK complexes.} It is still not known whether or not any of these  gives a triangulation of $\H P^2$. It has also been proved that these three combinatorial 8-manifolds are not isomorphic to each other though they have homeomorphic spaces.

\bigskip In this paper we will prove that none of these BK complexes can be embedded in the Euclidean space $E^{12}$ though all of them have a tight polyhedral embedding in $E^{14}$.  Since these BK complexes are simplicial complexes on 15 vertices, they are subcomplexes of the simplex $\Delta^{14}$ (a simplex with 15 vertices), it is obvious that each of these has a polyhedral embedding in $E^{14}$.

\section{Preliminaries}
 The three 8-dimensional combinatorial manifolds on $15$ vertices explicitly constructed in~\cite{bk1}, which we have called  BK complexes, are cohomology quaternionic projective planes. In fact, it has been proved there that any simplicial complex with 15 vertices, which is not a sphere, must be a cohomology quaternionic projective plane. It must be remarked here that the above three BK complexes are distinct because they have non-isomorphic automorphism groups: $M^{8}_{15}$ is characterized uniquely by the automorphism group $A_5$ acting transitively on the vertex set  and $\widetilde{M}^{8}_{15},\;\widetilde{\widetilde M^{8}_{15}}$ are characterized by the automorphism group $A_4$ and $S_5$ respectively. Hence all of these are combinatorially non-isomorphic. It is also proved there that these 8-manifolds are PL-homeomorphic to each other. As a consequence, therefore, the cohomology quaternionic projective plane $\H P^2$ does not have a unique triangulation on 15 vertices. This is in contrast with  real projective plane $\R P^2$ and the complex projective plane $\C P^2$ which have unique triangulations on 6 vertices and 9 vertices respectively.

  In this paper we follow some of the notations from \cite{mat}. The symbol [n] denotes the set $\{1,2,....,n\}$. For a simplicial complex $K$, the notation $(K*K)_\Delta$ denotes the deleted join of $K$ with itself. For the definition of deleted join of two simplicial complexes and of two topological spaces see (\cite{mat}, p.73). Cohomology groups are the simplicial cohomology groups with coefficients in $Z_2$. We use the $Z_2$-index of a $Z_2$-space as given in (\cite{mat}, p.95).

  Let $K$ be any simplicial complex on $n$ vertices $[n]=\{1,2,3\cdots , n\}$. The Alexander dual of $K$ is the simplicial complex $B(K)$ on the same vertex set that consists of all complements of all nonsimplices of  $K$, i.e.,
$$B(K):=\{G\subseteq \{1,2,\cdots ,n\}\;\overline{G}\notin K\}.$$
Here $\overline{G}$ denotes the complement of $G$ in the set $[n]$. The deleted join $(K\ast B(K))_{\Delta}$, denoted by $\hbox{Bier}_n(K)$, is called the Bier sphere on $n$ vertices $\{1,2,\cdots ,n\}$. Then we have the following result (see ~\cite{mat} Theorem 5.6.2 for a proof):

\begin{Thm}[Bier] The deleted join $\hbox{Bier}_n(K)=(K\ast B(K))_{\Delta}$ of $K$ is indeed a $(n-2)$ sphere with at most $2n$ vertices.\end{Thm}

\section{The complementarity of the BK complexes}

\bigskip\noindent Let $K$ denote any one of the 8-dimensional BK complexes each of which is a cohomology quaternionic projective plane $\H P^2$. In order to have an idea of the size of these complexes, let us assume that these are not 8-spheres. We will soon show that these are 5-neighborly, and so it is a consequence of Dehn-Sommerville equations that
$ f_0=\binom{15}{1}=15, f_1=\binom{15}{2}=105, f_2=\binom{15}{3}=455, f_3=\binom{15}{4}=1365,
f_4=\binom{15}{5}=3003, f_5=4515, f_6=4230, f_7=2205, f_8=490,$
where $f_i\; (i=0,1,\cdots ,8)$ denote the number of i-dimensional faces of $K$.

 Since Brehm and K\"{u}hnel, using the result of Marin \cite{mar}, have given only a sketch of the fact that $K$ satisfies the complementarity condition~\cite{bk1}, and since this property is crucial to our results in this paper, we give here a complete proof of the complementarity of $K$.

\bigskip Recall that a simplicial complex $K$ is said to satisfy  the {\it complementarity condition} if given any subset $X^0$ of the vertex set $K^0$, either $X^0$ is a simplex of $K$ or $K^0\setminus X^0$ is simplex of $K$ (mutually exclusive). Assuming that the Brehm-K\"{u}hnel complex $K$ is 5-neighborly (we prove it later), take any $X^0\subset K^0$ such that $|X^0|=5$. Since $|K^0\setminus X^0|=10$, and since $K$ cannot have a 9-simplex ($K$ is 8-dimensional), $K$ cannot have $K^0\setminus X^0$ as a simplex. Hence $X^0$ must be simplex of $K$. If $|X^0|<5$, then by the same argument $|K^0\setminus X^0|>10$ and $K^0\setminus X^0$ cannot be a simplex of $K$. It follows that complementarity for subsets of $K^0$ having at most 5 elements is now clear. However, if $|X^0|>5$, the argument is no longer valid. We now prove below that $K$ is not only 5-neighborly, but also that it satisfies the complementarity condition for any subset $X^0$ of $K^0$ whatever be the number of elements in $X^0$.

\bigskip\noindent Let us note the following:

\begin{Prop}Let $(K,h)$ be a cohomology quaternionic projective plane (mod 2) and $Y\subset K$ be a full subcomplex of $K$. Suppose $h|Y=0$ and $h_Y$ is the restriction of $h$ on the opposite complex  $K_Y$ in $K$. Then $(K_Y, h_Y)$ is also a cohomology quaternionic projective plane (mod 2).\end{Prop}

Here $(K,h)$ means $K$ is the BK complex on 15 vertices and $h\in H^4(K, \Z_2)$ is such that $h^2\neq 0$. The restriction is the homomorphism $H^4(K, \Z_2)\to H^4(K_Y, \Z_2)$ induced by the inclusion map $K_Y\subset K$.

\medskip

\noindent{\it Proof:}
 Note that $|Y|$ is a deformation retract of $U=|K|\setminus |K_Y|$ and $|K_Y|$ is a deformation retract of $V=|K|\setminus |Y|$. Hence $U,V$ are open sets and $|K|=U\cup V.$  If $h_Y=0,$  then $h/U=0$ and $h/V=0.$  But this means there exists an $\overline{h} \in H^2(K)$ such that $\overline{h}/U=h/U=0$ and $\overline{h}/V=h/V=0,$ which means $\overline{h}=0 \in H^2(|K|,U\cup V)$ and therefore $h\in H^2(|K|)=0$, a contradiction. \hfill $\Box$

\bigskip Let $L$ be the simplicial complex consisting of all simplices of $K$ and all 4-simplices formed by any five vertices from $K$. Then $L$ is a 5-neighborly simplicial complex which contains $K$ as a subcomplex.

\bigskip Let us define a 4-cochain $c$ on $L$ as follows: For any 4-simplex $\sigma$ of $L$, we put
$$c(\sigma)=\#\;\hbox{of 8-simplices of $L$ which are disjoint from $\sigma$(mod 2).}$$
Then extend $c$ to all 4-chains $X=\sum_{i}\sigma_i$ by $c(X)=\sum_i c(\sigma_i)$. This 4-cochain $c$ of $L$ is called a {\it counting 4-cochain} on $L$.

\bigskip\noindent {\bf Claim : } The counting 4-cochain $c$ is, in fact, a 4-cocycle.

\noindent {\it Proof : }
 We must prove that  if $\Delta$ is any 5-simplex in $L$, then $(\partial c)(\Delta)=c(\partial\Delta)=0$. Here $X=\partial\Delta$ is the boundary complex of $\Delta$ and hence is a 4-sphere. Note that there is no 8-simplex of $L$ disjoint from $\partial\Delta$ because if there is one, say $\sigma$, then $\sigma =L_\tau$ (the opposite complex) for some simplex $\tau$ in $\partial\Delta$. This will mean that $L_\tau ,\tau $ both are full subcomplexes of L and each one is contractible. But this will contradict Proposition 3.1 .\hfill $\Box$

\bigskip\noindent Now we have

 \begin{Prop} Suppose $\Sigma$ is a simple 4-sphere in a 5-neighborly simplicial complex $L$. Then $(c,\Sigma)\neq 0$ iff there is an 8-simplex $\Delta$ disjoint from $\Sigma$.\end{Prop}

\noindent {\it Proof : }
Suppose all the 8-simplices of $L$ intersect $\Sigma$. Note that such an 8-simplex $\Delta$ can have  only one vertex in common with $\Sigma$ because otherwise $\Delta$ will intersect all 4-simplices of $\Sigma$ and will not contribute anything to $(c,\Sigma)$. If $\Delta$ intersect in one vertex,  say $x$ of $\Sigma$, then $\Delta=x\ast\tau$ (join), where $\tau$ is a face of $\Delta$. This $\tau$ must be a face of an even number of such $\Delta$'s since the sum of all 8-simplices of $L$ forms a cycle. Thus for each $\sigma$ in $\Sigma$, there are even number of 8-simplices which intersect $\Sigma$ and so $(c,\sigma)=0$. Since the number of $\sigma$'s
in $\Sigma$ is even, we find that $(c,\Sigma)=0$.

 On the other hand, suppose there is an 8-simplex $\Delta$ which is disjoint from $\Sigma$. Then the opposite complex $L_{\Sigma}$ is the same as $\Delta$ and there must be only one $\Delta$ which is disjoint from $\Sigma$. This $\Delta$ does not contribute anything to $(c,\Sigma)$. Rest of all 8-simplices will intersect $\Sigma$ and there number is $\binom{9}{5}-1=125,$ an odd integer. Hence $(c,\Sigma)=1$.
\hfill $\Box$

 \begin{Lemma} The counting 4-cocycle $c$ has the property that the cup product $\overline{c}\cup\overline{c}\neq 0$ where $\overline{c}$ is the cohomology class of c.\end{Lemma}

\noindent {\it Proof : }
 Let $\Delta$ be an 8-simplex of $L$. Then $L_{\Delta}$ is a simple 4-sphere and $(h,L_{\Delta})=1$. But we also have that $(c,L_{\Delta})=1$. Hence $(h+\overline{c},L_{\Delta})=0$. This implies that $(h+\overline{c})\cup (h+\overline{c})=0$, which means $h\cup h+h\cup\overline{c}+\overline{c}\cup h+\overline{c}\cup\overline{c}=0$. Hence by anticommutativity of the cup product $h\cup h=\overline{c}\cup\overline{c}\neq 0$.
\hfill $\Box$

\begin{Prop} $K$ is 5-neighborly. \end{Prop}

\noindent {\it Proof : }
It suffices to prove that $L\subseteq K$. Let $\sigma$ be a 4-simplex in $L$. Then by definition of $L$ there is a simple 4-sphere $\Sigma$ containing $\sigma$. If $(c,\Sigma)=1$, then there is an 8-simplex $\Delta$ disjoint from $\Sigma$ in $L$ and hence in $K$. Hence by Proposition 3.1 applied to $K$ and $Y=\Delta$, the subcomplex $K_{\Delta}$ opposite to $\Delta$ in $K$ has a nontrivial $H^4$ and so $K_{\Delta}=\Sigma$, i.e., the 4-simplex $\sigma$ in $K$. On the other hand, if $(c,\Sigma)=0$, then apply Proposition 3.1 to $L$, $h=\overline{c}$ and let $Y$ be the full subcomplex $\overline{\Sigma}$ generated by $\Sigma$. Then we get that the restriction of $\overline{c}$ to $L_{\overline{\Sigma}}$ contains a simple 4-sphere $\Sigma'$ with $(c,\Sigma')=1$. Then Proposition 3.2 gives us an 8-simplex $\Delta$ disjoint from $\Sigma'$, it contains $\Sigma$ so that 4-simplex $\sigma$ is in $\Delta$ and is therefore in $L$. Thus $K=L$ and $K$ is 5-neighborly.
\hfill $\Box$

\begin{Prop} $K$ satisfies the complementarity condition.\end{Prop}

\noindent {\it Proof : }
 Let $X^0$ be any subset of $K^0$. Consider the full subcomplex generated by $X^0$ and by $K^0\setminus X^0$. By Proposition 3.1, one of them must have nontrivial $H^4$. Let it be $X^0$. Hence there is a nonbounding 4-sphere $\Sigma$ in $X^0$. Hence $K_{\Sigma}$ must be an 8-simplex whose vertex set contains $K^0\setminus X^0$. But this implies that $K^0\setminus X^0$ is a simplex in $K$.
\hfill $\Box$

\section{Main Result}

Now we prove the following:
{\Thm\label{main}  Let $K$ be any of the three combinatorial 8-dimensional BK complex each of which is a cohomology quaternionic projective plane $\H P^2$. Let $f:|K|\to E^{12}$ be a continuous map. Then there exists a pair $\sigma_1,\;\sigma_2$ of faces of $K$ such that
$$f(||\sigma_1||)\cap f(||\sigma_2||)\neq\phi$$}

{\Cor \rm None of the 8-dimensional BK complexes can be embedded in $E^{12}$.
}

\bigskip
\noindent {\it Proof of the Theorem 4.1:}
 Suppose for every disjoint pair $\sigma_1,\;\sigma_2$ of faces of $K,\;f(||\sigma_1||)\cap f(||\sigma_2||) = \phi$. Let $K^{\ast 2}_{\Delta}$ denote the deleted join $(K\ast K)_{\Delta}$ of $K$ and $({E^{12}})^{\ast 2}_{\Delta}$ denote the deleted join of $E^{12}$. Then clearly both are free $\Z_2$-spaces (interchange of the coordinates) and
$$f\ast f :K^{\ast 2}_{\Delta} \to ({E^{12}})^{\ast 2}_{\Delta}$$ defines a $\Z_2$-map (see~\cite{mat} p.108 ). This means $$\hbox{ind}_{\Z_2} K^{\ast 2}_{\Delta}\leq \hbox{ind}_{\Z_2}({E^{12}})^{\ast 2}_{\Delta}.\hspace{2in} (\ast)$$ It is well-known (see~\cite{mat} p.110 Lemma 5.5.4) that $\hbox{ind}_{\Z_2}({E^{12}})^{\ast 2}_{\Delta}=12$. On the other hand we now compute the $\hbox{ind}_{\Z_2} K^{\ast 2}_{\Delta}$.

 Note that by Proposition 3.4 the simplicial complex $K$ satisfies the complementarity property. Let us consider the Bier sphere $B_{15}(K)$. Now if $F\in K$, then its complement $K^0\setminus F^0$ is not a simplex of $K$ and hence by the definition of $B_{15}(K),\;F\in B_{15}(K)$.
Conversely if $F\in B_{15}(K)$, then $K^0\setminus F^0$ is a nonsimplex of $K$ and hence by complementarity, $F^0$ must form a  simplex of $K$, i.e., $F\in K$ iff $F\in B_{15}(K)$. In other words $K=B_{15}(K)$.

It now follows from Bier's Theorem that $(K\ast B_{15}(K))_{\Delta}=K^{\ast 2}_{\Delta}$ is a 13-sphere and hence $\hbox{ind}_{\Z_2} K^{\ast 2}_{\Delta}=13$. But this contradicts $(\ast)$ and proves the Theorem. \hfill $\Box$

\bigskip
We now indicate yet another proof of our result. The following theorem is well-known (see \cite{mat} p.120)

\begin{Thm}[Sarkaria's coloring/embedding theorem] Let $K$ be simplicial complex on n vertices and let $ \mathfrak{F}= \mathfrak{F}(K)$  be the system of minimal nonfaces of $K$. Then

$$ ind_{Z_2}(K^{*2}_{\Delta} ) \geq n - \chi (KG(\mathfrak{F}))-1,  $$
where $KG(\mathfrak{F})$ denotes the Kneser graph of the family $\mathfrak{F}$ and $\chi (KG(\mathfrak{F}))$ denotes the chromatic number of the Kneser graph (see ~\cite{mat}, p.119, for definitions of a set system and the Kneser graph of a set system).
\end{Thm}

\noindent{\it Proof of Theorem 4.1 : }  In our case each of the BK complexes $K$  is a simplicial complex on 15 vertices. The nonfaces are the complements of all 5-simplexes, all 6-simplexes, all 7 simplexes and all 8-simplexes of $K$.These nonfaces are also minimal because if $G$ is such a nonface and $G\cup A$ is a face of $K$ for any nonempty subset $A$ of [n], then $G$ must be a simplex of $K$, a contradiction.

\medskip
\noindent{\bf Claim:} If $G_1, G_2 \in \mathfrak{F}$, then $G_1\cap G_2 \neq \phi .$

\medskip
\noindent This follows from the complementarity property of $K$. Since $G_2$ is a nonface of $K$,  $F=[n]\setminus G_2$ is a face of $K$. Hence if $G_1\cap G_2 =\phi$, then $G_1\subset F$ is also a face of $K$, a contradiction. It follows that no two members of $\mathfrak{F}$  are disjoint, i.e., the Kneser graph $KG(\mathfrak{F})$ does not have any edge and hence the chromatic number $\chi(KG(\mathfrak{F})) =1.$ Therefore by the Sarkaria's Theorem quoted above K cannot be embedded in the Euclidean space $E^d$  for $d\leq n- \chi(KG(\mathfrak{F}))-2= 12.$  \hfill $\Box$

{\Remark \rm The well known J. Radon's Theorem says that if there is a continuous map $f:\partial\Delta^{n+1}=K\to E^n\;(n=1,2,\cdots )$, then there is a pair of disjoint faces $\sigma_1,\sigma_2$ of $K$ such that $f(||\sigma_1||)\cap f(||\sigma_2||)\neq\phi$. Our theorem can be seen as an interesting analogue of Radon's theorem from $\H P^1=\S^4$ to three $BK$ complexes which are cohomology quaternionic projective planes.
}

{\Remark \rm A method similar to our Theorem~\ref{main} can also be applied to prove the nonembedding of the cohomology Cayley plane, but at present it is not known whether or not there is a 16-dimensional cohomology Cayley plane on 27 vertices.}

\end{document}